\documentclass[12pt]{amsart}
\usepackage{amssymb}

  \theoremstyle{definition}
  \newtheorem{theorem}{Theorem}[section]
  \newtheorem{corollary}[theorem]{Corollary}

  \newtheorem{example}[theorem]{Example}

  \theoremstyle{remark}
  \numberwithin{equation}{section}



\begin{document}
\bibliographystyle{plain}
\title[Subspaces of $L_p$ with unconditional bases]{Subspaces of $L_p$, $p>2$, with unconditional bases have equivalent partition and
weight norms}
   
\vskip
2 cm
\author {Dale E. Alspach}
\address{Department of Mathematics, Oklahoma State
   University, Stillwater, OK 74078}
\email{alspach@math.okstate.edu}
\author {Simei Tong}
\address{Department of Mathematics, University of Wisconsin-Eau
   Claire, Eau Claire, WI 54702}
\email{tongs@uwec.edu}

\begin{abstract} In this note we give a simple proof that every subspace
of $L_p$, $2<p<\infty$, with an unconditional basis has an equivalent norm
determined by
partitions and weights. Consequently $L_p$ has a norm determined by partitions
and weights.
\end{abstract}
\maketitle

\section{Introduction}
It is a well known consequence of Khintchine's inequality
that if $(x_n)$ is an unconditional basic sequence in $L_p=L_p[0,1]$
with sign unconditional constant $K$, then
$$
\frac{1}{A_pK}\left ( \int \left (\sum_n a_n^2 x_n^2 \right )^{p/2}\right )^{1/p} 
\le \|\sum_n a_n x_n\|_p \le B_p K\left ( \int \left (\sum_n a_n^2 x_n^2 \right )^{p/2}\right )^{1/p},
		$$
		where $A_p$ and $B_p$ are the Khintchine constants, \cite[Inequality 1.7]{Lp}. In this note we show that this inequality easily implies the following.

		\begin{theorem} Suppose that $X$ is a subspace of $L_p$, $2<p<\infty,$ with
		an unconditional basis. Then there is an equivalent norm on $X$ given
by partitions and weights.
\end{theorem}

and

\begin{corollary} $L_p[0,1]$, $2<p<\infty$,
 has an equivalent norm given by partitions and weights.
\end{corollary}

These results answer some questions raised by the authors, \cite{AT}. In
that paper and the second author's dissertation, \cite{Tong}, the notion of norm
given by partition and weights was introduced and some development of
the concept was made.

Define for any partition
$P=\{ N_i\}$ of $\mathbb{N}$ and weight function $ W : {\mathbb{N}}
 \rightarrow { ( 0, 1 ]}$
  $$ \|(a_i)\|_{P,W} = \left(\sum_{i} \left (
 \sum_{j\in {N_i}} {a_j}^2 w_{j}^{2}\right)^{\frac{p}{2}}\right)
^{\frac{1}{p}} $$
Now suppose that ${\mathcal{P}}=(P_k, W_k)_{k\in K}$ is a family of
pairs of
partitions and functions as above.
Define
$$
 \|(a_i)\|_{\mathcal{P}} =  \sup_{k\in K}
 \|(a_i)\|_{(P_k, W_k)}.
$$
We say that
$\|(a_i)\|_{\mathcal{P}}$ is a norm determined by partitions and weights.

\section{Proof of the Theorem}

\begin{proof} Let $X$ be a subspace of $L_p$ with normalized
unconditional basis $(x_n)$. Let $G$ be the set of all $g \in L_{p/(p-2)}$,
$g$ non-negative, and $\|g\|_{p/(p-2)}=1.$ For each $g \in G$
 let $w_{g,n}=(\int g x_n^2)^{1/2}$ for $n =1, 2, \dots.$
Now suppose that $(a_n)$ is a sequence of real numbers with only finitely many non-zero.
Then
\begin{eqnarray*}
\| \sum a_n x_n\| &\sim &
\left ( \int \left (\sum_n a_n^2 x_n^2 \right )^{p/2}\right )^{1/p} \\
&= &\|\sum a_n^2 x_n^2 \|_{p/2}^{1/2} \\
&= &\left (\sup\{\int g\sum a_n^2 x_n^2:g \in G\} \right )^{1/2}\\
&= &\left (\sup \{\sum a_n^2 w_{g,n}^2: g\in G\} \right )^{1/2}.\\
&= &\sup \{ \left (\sum a_n^2 w_{g,n}^2  \right )^{1/2} : g\in G\}
\end{eqnarray*}

It follows immediately that if $\mathcal P=\{(P_I,(w_{g,n})):g \in G\},$ where
$P_I$ denotes the trivial partition,
then 
$$ \||\sum a_n x_n\|| = \|(a_n)\|_{\mathcal P}$$
is an equivalent norm on $X$.
\end{proof}

Note that in the proof we really have a correspondence between norm one
positive operators from $L_{p/(p-2)}$ into $\ell_\infty$
and unconditional basic sequences in $L_p$. Indeed, if $T$ is such an
operator, then for each $g \in L_{p/(p-2)}$, $0<\|g\|\le 1,$
and $g \ge 0$ let $(w_{g,n})_{n=1}^\infty=(Tg(n)^{1/2})$. Let $(e_n)$ be
the coordinate functionals on $\ell_\infty$ and $x_n = T^*(e_n)^{1/2}$
for all $n$. (We assume $x_n \ne 0.$)
Then $(x_n \otimes r_n)$, where $(r_n)$ is the sequence of
Rademacher
functions, is an unconditional basic sequence (not necessarily normalized)
in $L_p([0,1]^2)$ with equivalent norm given by the family of partitions and
weights $(P_I, (w_{g,n}))_{g\in G}$.

In the proof we used the entire
set $G$, but it easy to see that it would be sufficient to use a norm
dense subset of $G$. Also the proof yields only the indiscrete partition.
In some of the examples below we show how other partitions arise naturally.

Now we consider some examples.

\begin{example}
Suppose $x_n=\lambda(A_n)^{-1/p}1_{A_n}$ for all $n \in \mathbb N$
and $A_m \cap A_n = \emptyset$ if $n \ne m$. If $g \in L_{p/(p-2)}$,
then $\int g x_n^2 =\int E(g|\mathcal F)x_n^2$
where $E(g|\mathcal F)$ is the conditional expectation with respect to the
$\sigma$-algebra $\mathcal F$
generated by the sets $(A_n)$. Thus we may assume $g =
\sum b_k \lambda(A_k)^{-(p-2)/p} 1_{A_k}$ for some scalars $(b_k)$. Then

$$\int g x_n^2 = b_n \lambda(A_n)^{-2/p-(p-2)/p+1}\\ =b_n .$$
 Therefore
$w_{g,n}= b_n^{1/2}$ for all $n$ and
$$\||\sum a_n x_n\||^2= \sup\{\sum
a_n^2 b_n :g \in G\}= \sup \{\sum a_n^2 b_n: \sum b_k^{p/(p-2)}=1,b_k
\ge 0\} = (\sum a_n^{p})^{2/p}.$$
Thus we get $\ell_p$ isometrically.
\end{example}

The previous example easily generalizes to any sequence of disjointly supported functions.

\begin{example}
Let $(A_n)$ be a sequence of disjoint sets of positive measure and for each $n$ let $x_n$ be a norm one element suported in $A_n$. Arguing similarly to the previous example (or by change of measure) we may assume that $g = \sum_n b_k y_k$ where  for each $k$, $y_k$ is the unique norm one function in $L_{p/(p-2)}$ which is norming for $x_k^2.$ Then
$\int g x_n^2 = b_n$ and we conclude as before that we have an isometric copy of $\ell_p.$
\end{example}

\begin{example} Suppose $x_{n,j}=\lambda(A_n)^{-1/p}1_{A_n}r_j$ for any $n
\in \mathbb N$, $j \in \mathbb N$, $(r_j)$ is the sequence of Rademacher
functions and $A_m \cap A_n = \emptyset$ if $n \ne m$. Because $x_{n,j}^2
=x_n^2$ in the first example, it follows that

\begin{multline*}
\||\sum_n\sum_j a_{n,j} x_{n,j}\||^2=\sup \{\sum_n \sum_j a_{n,j}^2
w_{g,n,j}^2 :g\in G\}= \sup\{\sum_n \sum_j a_{n,j}^2 b_n :g \in G\}\\
= \sup \{\sum_n \sum_j a_{n,j}^2 b_n: \sum b_k^{p/(p-2)}=1,b_k \ge 0\} =
(\sum_n(\sum_j a_{n,j}^2)^{p/2})^{2/p}.  \end{multline*}

Thus we get $(\sum \ell_2)_{\ell_p}$ isometrically.
\end{example}

By arguing as in the second example it follows that for an unconditional basic sequence
which breaks into a sequence of disjointly supported (on $[0,1]$)
unconditional basic subsequences the norm $\||\cdot \||$ will be an isometric $\ell_p$ sum of the norms for the subsequences.

\begin{example}
Let $(x_n)$ 
be a sequence of stochastically independent normalized functions in $L_p$.
Now suppose that $g \in G$ and $(a_n)$ is a sequence of real numbers
with only finitely many non-zero. Then by an inequality of Rosenthal
\cite[Lemma 1]{MR42:6602}
\begin{multline*} \sup_{g \in G} \sum_n
a_n^2 \int g x_n^2 =\sup_{g \in G} \int g \sum a_n^2 x_n^2 = \|\sum
a_n^2 x_n^2 \|_{p/2} \\ \le K \max \{\|\sum a_n^2 x_n^2 \|_1, \left (
\sum \|a_n^2 x_n^2\|_{p/2}^{p/2}\right )^{2/p}\} = K\max \{ \sum a_n^2
\|x_n\|_2^2, \left (\sum |a_n|^p \|x_n\|_p^p\right )^{2/p}\}.  \end{multline*}
Hence $\||\sum a_n x_n \|| \le K\max \{ \left (\sum a_n^2 \|x_n\|_2^2\right)^{1/2}, \left (\sum |a_n|^p\right )^{1/p}\}.$

Suppose $(c_n) \in \ell_{p/(p-2)}$ with $\|(c_n)\|=1$, $c_n \ge 0$ for all $n$, $c_n$ non-zero for only finitely many $n$. Let
$g = \max c_n |x_n|^{p-2}.$ Because $g=\sum x_n {}_{|\{\tau=n\}}$,
where $\tau(t) = \min \{k:c_k |x_k|^{p-2}(t)=g(t)\}$, $\|g\|_{p/(p-2)}\le (\sum
c^{p/(p-2)}_n \|x_n\|_p^p)^{(p-2)/p}=1.$ Then $\int g x_n^2 \ge \int c_n x_n^p
= c_n.$ Therefore $\||\sum a_n x_n \|| \ge \sup_{\|(c_n)\|_{p/(p-2)}=1}
(\sum a_n^2 c_n)^{1/2} = (\sum |a_n|^p)^{1/p}.$ If $g = 1$, then
$w_{g,n}^2=\int x_n^2 =\|x_n \|_2^2.$  Thus \newline $\||\sum a_n x_n \||
 \ge \max \{ \left (\sum a_n^2 \|x_n\|_2^2\right)^{1/2}, \left (\sum
 |a_n|^p\right )^{1/p}\}.$
\end{example}

The estimate using $g=\max c_n |x_n|^{p-2}$
above actually gives us some information about
any normalized unconditional basis $(x_n)$ in $L_p,$
$$\sup\{(\sum a_n^2 w_{g,n}^2)^{1/2}: g=\max c_n |x_n|^{p-2},\|(c_n)\|_{p/(p-2)}=1,
(c_n)\ge 0\} \ge (\sum |a_n|^p)^{1/p}.$$ Thus we can always include the
discrete partition $P_D=\{\{n\}:n\in \mathbb N\}$ with constant weight one
in the family of partitions and weights $\mathcal P.$

Finally we consider the Haar system. 
\begin{example}
Let
$h_{n,k} = 2^{(n-1)/p}1_{[2^{-n}2k,2^{-n}(2k+1))}-
2^{(n-1)/p}1_{[2^{-n}(2k+1),2^{-n}(2k+2))}$ for $n= 1,2, \dots $, $k = 0,
1, \dots 2^{n-1}-1.$ Let $h_{0,0} =1_{[0,1)}.$

Suppose $g \in G$, $n\in \mathbb N$, and $g =
\sum_{k=0}^{2^{n-1}-1}b_{n,k}2^{(n-1)((p-2)/p)}1_{[2^{1-n}k,2^{1-n}(k+1))}.$ Then
for $0<m \le n$,
$$\int g h_{m,\ell}^2 = \sum_{k=\ell2^{n-m}}^{(\ell+1)2^{n-m}-1}
\int
b_{n,k}2^{(n-1)(p-2)/p+2(m-1)/p}1_{[2^{1-n}k,2^{1-n}(k+1))}
=\sum_{k=\ell 2^{n-m}}^{(\ell+1)2^{n-m}-1} b_{n,k}2^{2(m-n)/p}.$$
For $m=0$, $\int g h_{0,0}^2= \int g h_{1,0}^2=\sum_{k=0
}^{2^{n-1}-1} b_{n,k}2^{2(1-n)/p}.$
For $m>n$ and \newline
$k 2^{m-n}\le \ell <(k+1)2^{m-n}$,
$$\int g h_{m,\ell}^2 = \int
b_{n,k}2^{(n-1)(p-2)/p+2(m-1)/p}1_{[2^{1-m}k,2^{1-m}(k+1))}
= b_{n,k} 2^{(n-m)(p-2)/p}.$$
Thus
\begin{multline*} \|| a_{0,0}h_{0,0}+\sum_{m=1}^\infty \sum_{\ell=0}^{2^{m-1}-1}
a_{m,\ell} h_{m,\ell}\||
=\sup\left (
a_{0,0}^2\sum_{k=0
}^{2^{n-1}-1} b_{n,k}2^{2(1-n)/p}\right . \\
\left .+\sum_{m=1}^n \sum_{\ell=0}^{2^{m-1}-1}
a_{m,\ell}^2 \sum_{k=\ell
2^{n-m}}^{(\ell+1)2^{n-m}-1} b_{n,k}2^{2(m-n)/p}+\sum_{m=n+1}^\infty
\sum_{\ell=0}^{2^{m-1}-1}
a_{m,\ell}^2 b_{n,k} 2^{(n-m)(p-2)/p}\right )^{1/2}
\end{multline*}
where the supremum is over all $n$ and all non-negative
sequences
$(b_{n,k})_{k=0}^{2^{n-1}-1}$\newline
 with $\|(b_{n,k})\|_{(p-2)/p}
=1.$

The contribution from the terms $m>n$ can be neglected,
so an equivalent norm on $L_p$ is 
$$\sup \left (a_{0,0}^2\sum_{k=0
}^{2^{n-1}-1} b_{n,k}2^{2(1-n)/p}+ \sum_{m=1}^n \sum_{\ell=0}^{2^{m-1}-1} a_{m,\ell}^2
2^{2(m-n)/p}\sum_{k=\ell 2^{n-m}}^{(\ell+1)2^{n-m}-1} b_{n,k}\right
)^{1/2}$$
where the supremum is over the same sequences as above.
\end{example}

\bibliography{partw}
\end{document}